\documentclass[leqno,12pt]{article}
\usepackage{amssymb,amsfonts}
\usepackage{amsmath,latexsym}

\usepackage{amstext,epic,eepic,epsf,pslatex}
\usepackage{graphicx}

\newcommand{\bepr}{{\em Proof} } 
\newcommand{\enpr}{\hfill \rule{.5em}{.5em}}

\newcommand{\R}{{\mathbb R}}
\newcommand{\N}{{\mathbb N}}

\def\XXint#1#2#3{{\setbox0=\hbox{$#1{#2#3}{\int}$ }
\vcenter{\hbox{$#2#3$ }}\kern-.6\wd0}}

\newtheorem{thm}{Theorem}[section] 
\newtheorem{lemma}{Lemma}[section]

\newtheorem{cor}{Corollary}[section]

\begin{document}

\title{Multi-dimensional Burgers equation with unbounded initial data: well-posedness and dispersive estimates}

\author{Denis Serre\thanks{U.M.P.A., UMR CNRS--ENSL \# 5669. \'Ecole Normale Sup\'erieure de Lyon. {\tt denis.serre@ens-lyon.fr}}
\\
Luis Silvestre\thanks{Department of Mathematics, the University of Chicago. {\tt luis@math.uchicago.edu}}
\thanks{LS was partially supported by NSF grants DMS-1254332 and DMS-1764285.}}

\date{}

\maketitle

\begin{abstract}
The Cauchy problem for a scalar conservation laws admits a unique entropy solution when the data $u_0$ is a bounded measurable function (Kruzhkov). The semi-group $(S_t)_{t\ge0}$ is contracting in the $L^1$-distance.

For the multi-dimensional Burgers equation, we show that $(S_t)_{t\ge0}$ extends uniquely as a continuous semi-group over $L^p(\R^n)$ whenever $1\le p<\infty$, and $u(t):=S_tu_0$ is actually an entropy solution to the Cauchy problem. When $p\le q\le \infty$ and $t>0$, $S_t$ actually maps $L^p(\R^n)$ into $L^q(\R^n)$.

These results are based upon new dispersive estimates. The ingredients are on the one hand Compensated Integrability, and on the other hand a De Giorgi-type iteration.
\end{abstract}



\paragraph{Notations.} When $1\le p\le\infty$, the natural norm in $L^p(\R^n)$ is denoted $\|\cdot\|_p$, and the conjugate exponent of $p$ is $p'$. The total space-time dimension is $d=1+n$ and the coordinates are $x=(t,y)$. In the space of test functions, ${\cal D}^+(\R^{1+n})$ is the cone of functions which take non-negative values.  The partial derivative with respect to the coordinate $y_j$ is $\partial_j$, while the time derivative is $\partial_t$. Various finite positive constants that depend only the dimension, but not upon the solutions of our PDE, are denoted $c_d,c_{d,p}, c_{d,p,q}$~; they usually differ from one inequality to another one. We denote $C_0(0,+\infty)$ the space of continuous functions over $(0,+\infty)$ which tend to zero at infinity. Mind that $C(\R_+)$ is the space of bounded continuous functions over $[0,+\infty)$.

\section{Introduction}

Let us consider a scalar conservation law in $1+n$ dimensions
\begin{equation}
\label{eq:mdcl}
\partial_tu+\sum_{i=1}^n\partial_if_i(u)=0,\qquad t>0,\,y\in\R^n.
\end{equation}
We complement this equation with an initial data
$$u(0,y)=u_0(y),\qquad y\in\R^n.$$
The flux $f(s)=(f_1(s),\ldots,f_n(s))$ is  a smooth vector-valued function of $s\in\R$. We recall the terminology that an entropy-entropy flux pair is a couple $(\eta,q)$ where $s\mapsto\eta(s)$ is a numerical function, $s\mapsto q(s)$ a vector-valued function, such that $q'(s)\equiv\eta'(s)f'(s)$. The Kruzhkov's entropies and their fluxes form a one-parameter family:
$$\eta_a(s)=|s-a|,\qquad q_a(s)={\rm sgn}(u-a)\,(f(u)-f(a)).$$
Together with the affine functions, they span the cone of convex functions. 

\bigskip

We recall that an {\em entropy solution} is a measurable function $u\in L^1_{\rm loc}([0,+\infty)\times\R^n)$ such that $f(u)\in L^1_{\rm loc}([0,+\infty)\times\R^n)$, which satisfies the Cauchy problem in the distributional sense,
\begin{equation}
\label{eq:distr}
\int_0^\infty dt\int_{\R^n}(u\partial_t\phi+f(u)\cdot\nabla_y\phi)\,dy+\int_{\R^n}u_0(y)\phi(0,y)\,dy=0,\qquad\forall\phi\in{\cal D}(\R^{1+n}),
\end{equation}
together with the entropy inequalities
\begin{eqnarray}
\nonumber
\int_0^\infty dt\int_{\R^n}(\eta_a(u)\partial_t\phi+q_a(u)\cdot\nabla_y\phi)\,dy & & \\
\label{eq:entrina}
+\int_{\R^n}\eta_a(u_0(y))\phi(0,y)\,dy & \ge0, & \qquad\forall\,\phi\in{\cal D}^+(\R^{1+n}),\,\forall\,a\in\R.
\end{eqnarray}

\bigskip

The theory of this Cauchy problem dates back to 1970, when S. Kruzhkov \cite{Kru} proved that if $u_0\in L^\infty(\R^n)$, then there exists one and only one entropy solution in the class
$$L^\infty(\R_+\times\R^n)\cap C(\R_+;L^1_{\rm loc}(\R^n)).$$
The parametrized family of operators $S_t:u_0\mapsto u(t,\cdot)$, which map $L^\infty(\R^n)$ into itself, form a semi-group. We warn the reader that $S_t:L^\infty\rightarrow L^\infty$ is not continuous, because of the onset of shock waves. Likewise, $t\mapsto u(t)$ is not continuous from $\R_+$ into $L^\infty(\R^n)$.

This semi-group enjoys nevertheless nice properties. On the one hand, a comparison principle says that if $u_0\le v_0$, then $S_tu_0\le S_tv_0$. For instance, the solution $u$ associated with the data $u_0$ is majorized by the solution $\bar u$ associated with the data $(u_0)_+$, the positive part of $u_0$. On another hand, if $v_0-u_0$ is integrable over $\R^n$, then $S_tv_0-S_tu_0$ is integrable too, and
\begin{equation}
\label{eq:contr}
\int_{\R^n}|S_tv_0-S_tu_0|(y)\,dy\le\int_{\R^n}|v_0-u_0|(y)\,dy.
\end{equation}
Finally, $S_t$ maps $L^p\cap L^\infty(\R^n)$ into itself, and the function $t\mapsto \|S_tu_0\|_p$ is non-increasing.

\bigskip

Because of (\ref{eq:contr}) and the density of $L^1\cap L^\infty(\R^n)$ in $L^1(\R^n)$, the family $(S_t)_{t\ge0}$ extends in a unique way as a continuous semi-group of contractions over $L^1(\R^n)$, still denoted $(S_t)_{t\ge0}$. When $u_0\in L^1(\R^n)$ is unbounded, we are thus tempted to declare that $u(t,y):=(S_tu_0)(y)$ is the {\em abstract solution} of the Cauchy problem for (\ref{eq:mdcl}) with initial data $u_0$. At this stage, it is unclear whether $(S_t)_{t\ge0}$ can be defined as a semi-group over some $L^p$-space for $p\in(1,\infty)$, because  the contraction property (\ref{eq:contr}) occurs only in the $L^1$-distance, but in no other $L^p$-distance.

An alternate construction of $(S_t)_{t\ge0}$ over $L^1(\R^n)$, based upon the Generation Theorem for nonlinear semigroups, was done by M. Crandall \cite{Cra}, who pointed out that  it is unclear  whether $u$ is an entropy solution, because the local integrability of the flux $f(u)$ is not guaranted\footnote{Except of course in the case where $f$ is globally Lipschitz.}. The following question is therefore an important one:
\begin{quote}
{\em  Identify the widest class of integrable initial data for which $u$ is actually an entropy solution of (\ref{eq:mdcl}). }
\end{quote}

\bigskip

Our most complete results are about a special case, the so-called {\em multi-dimensional Burgers equation}
\begin{equation}
\label{eq:Bru}
\partial_tu+\partial_j\frac{u^2}2\,+\cdots+\,\partial_n\frac{u^{n+1}}{n+1}=0,
\end{equation}
which is a paradigm of a genuinely non-linear conservation law. This equation was already considered by G. Crippa et al. \cite{COW}, and more recently by L. Silvestre \cite{Sil}.
The particular flux in (\ref{eq:Bru}) is a prototype for genuinely nonlinear conservation laws, those which satisfy the assumption
\begin{equation}
\label{eq:GNL}
\det(f'',\ldots,f^{(n+1)})\ne0.
\end{equation}
The latter condition is a variant of the {\em non-degeneracy condition} at work in the kinetic formulation of the equation (\ref{eq:mdcl})~; see \cite{LPT} or \cite{Per}.

Our first result deals with dispersive estimates:
\begin{thm}\label{th:Bcomp}
Let $1\le p\le q\le\infty$ be two exponents. Define two parameters $\alpha,\beta(p,q)$ by
\begin{equation}
\label{eq:aldef}
\alpha(p,q)=\frac{h(q)}{h(p)}\,,\qquad h(p):=2+\frac{dn}p\,
\end{equation}
and
\begin{equation}
\label{eq:bedef}
\beta(p,q)=h(q)(\delta(p)-\delta(q)),\qquad\delta(p):=\frac n{2p+dn}\,.
\end{equation}
There exists a finite constant $c_{d,p,q}$ such that
for every initial data $u_0\in L^1\cap L^\infty(\R^n)$, the entropy solution  $u(t)$ of the scalar conservation law (\ref{eq:Bru}) satisfies
\begin{equation}\label{eq:dispq}
\|u(t)\|_q\le c_{d,p,q}\,t^{-\beta(p,q)}\|u_0\|_p^{\alpha(p,q)},\qquad \forall \,t>0.
\end{equation}
\end{thm}

\bigskip

\paragraph{Remarks}
\begin{itemize}
\item The consistency of estimates (\ref{eq:dispq}) with the H\"older inequality is guaranted by the property that whenever $\theta\in(0,1)$,
\begin{equation}
\label{eq:albeH}
\left(\frac1q=\frac{1-\theta}p+\frac\theta{r}\right)\Longrightarrow
\left\{\begin{array}{l}
\alpha(p,q)=1-\theta+\theta\alpha(p,r), \\ \\ \beta(p,q)=\theta\beta(p,r).
\end{array}\right.
\end{equation}
\item The consistency under composition $(p,q)\wedge(q,r)\mapsto(p,r)$ is ensured by the rules
\begin{equation}
\label{eq:albecomp}
\alpha(p,r)=\alpha(p,q)\alpha(q,r)\qquad\hbox{and}\qquad
\beta(p,r)=\beta(q,r)+\beta(p,q)\alpha(q,r)
\end{equation}
\item In one space dimension, (\ref{eq:dispq}) gives back well-know results, such as Theorem\footnote{ Mind that this statement contains a typo, as the choice $r=1-\frac1p$ in Theorem 11.5.1 yields the exponent $-\frac1{p+1}$ instead of $-\frac{p}{p+1}$\,.} 11.5.2 in \cite{Daf_book}.
\end{itemize}

\bigskip

Theorem \ref{th:Bcomp} has several important consequences. An obvious one is that the extension of $(S_t)_{t\ge0}$ as a semi-group over $L^1(\R^n)$ satisfies the above estimates with $p=1$~:
\begin{cor}\label{c:dispun}
If $u_0\in L^1(\R^n)$ and $t>0$, then $S_tu_0\in\bigcap_{1\le q\le\infty} L^q(\R^n)$ and we have
$$\|S_tu_0\|_q\le c_{d,q}\,t^{-\kappa/q'}\|u_0\|_1^{1-\nu/q'},\qquad\forall q\in[1,\infty],$$
where the exponents are given in terms of
$$\kappa=2\,\frac{d-1}{d^2-d+2}\,\quad\hbox{and}\quad\nu=\frac{d(d-1)}{d^2-d+2}\,.$$
\end{cor}

The next one is that the Cauchy problem is solvable for data taken in $L^p(\R^n)$ for arbitrary exponent $p\in[1,\infty]$. In particular, it solves Crandall's concern.
\begin{thm}\label{th:wpp}
Let $p\in[1,\infty)$ be given. For every $t\ge0$, the operator $S_t:L^1\cap L^\infty(\R^n)\rightarrow L^1\cap L^\infty(\R^n)$ admits a unique continuous extension $S_t:L^p(\R^n)\rightarrow L^p(\R^n)$. 

The family $(S_t)_{t\ge0}$ is a continuous semi-group over $L^p(\R^n)$. If $u_0\in L^p(\R^n)$, the function $u(t,y)$ defined by $u(t)=S_tu_0$ is actually an entropy solution of the Cauchy problem for (\ref{eq:Bru}) with initial data $u_0$.

Finally, $S_t(L^p(\R^n))$ is contained in $\bigcap_{p\le q\le\infty}L^q(\R^n)$ and the estimates (\ref{eq:dispq}) are valid for every data $u_0$ in $L^p(\R^n)$.
\end{thm}

The proof of Theorem \ref{th:Bcomp} will be done in two steps. The first one consists in establishing the estimate (\ref{eq:dispq}) when $q=p^*$ is given by the formula
$$p^*=d\left(1+\frac pn\right).$$
To this end, we apply Compensated Integrability to a suitable symmetric tensor, whose row-wise divergence is a bounded measure with controlled mass. This argument involves the theory recently developped by the first author in \cite{Ser_DPT,Ser_CI}.
The second step is an iteration in De Giorgi's style, based on the preliminary work \cite{Sil} by the second author~; see also the original paper by E. De Giorgi \cite{DG} or the review paper by A. Vasseur \cite{Vas}. This technique allows us to establish an $L^\infty$-estimate, which extends the dispersive estimate to $q=+\infty$. Then using the H\"older inequality, we may interpolate between this result and the decay of $t\mapsto\|u(t)\|_p$, and treat every exponent $q>p$.

We notice that the symmetric tensor mentionned above extends to a multi-dimensional context the one already used when $n=1$ by L. Tartar \cite{Tar_HW} to prove the compactness of the semi-group, and by F. Golse \cite{Gol} (see also \cite{GoPe}) to prove some kind of regularity.

\bigskip

\paragraph{Previous dispersive estimates.}
In one space dimension $n=1$, (\ref{eq:Bru}) reduces to the original Burgers equation. Its Kruzhkov solution satisfies the Oleinik inequality $\partial_yu\le\frac1t$\,, which does not involve the initial data at all. Ph. B\'enilan \& M. Crandall \cite{BC}  proved
\begin{equation}
\label{eq:Daf}
TV\left(\frac{u(t)^2}2\right)\le\frac{2\|u_0\|_1}t\,,
\end{equation}
by exploiting the homogeneity of the flux $f(s)=\frac{s^2}2$\,. Inequality (\ref{eq:Daf}) implies an estimate
\begin{equation}
\label{eq:heat}
\|u(t)\|_\infty\le2\,\sqrt{\frac{2\|u_0\|_1}t\,}\,,
\end{equation}
which is a particular case of Corollary \ref{c:dispun} in this simplest case.

C. Dafermos \cite{Daf} proved a general form of (\ref{eq:Daf}) in situations where the flux $f$ may have one inflexion point and the data $u_0$ has bounded variations, by a clever use of the generalized backward characteristics.  His argument involves the order structure of the real line. Backward characteristics are not unique in general. Given a base point $(x^*,t^*)$ in the upper half-plane, one has to define and analyse the minimal and the maximal ones. The description of backward characteristics seems to be much more complicated in higher space dimensions, and Dafermos' strategy has not been applied successfully beyond the 1-D case.

\paragraph{Enhanced decay.} Because of a scaling property which will be described in the next section, the dispersion (\ref{eq:dispq}) is optimal, as long as we involve only the $L^p$-norms, and we exclude any extra information about the initial data. It is however easy to obtain a better decay as time $t$ goes to infinity. Let us give one example, by taking an initial data $u_0$ such that
$$0\le u_0(y)\le v_0(y_1),\qquad v_0\in L^1(\R).$$
By the maximum principle, we have $u(t,y)\le v(t,y_1)$, where $v$ is the solution of the $1$-dimensional Burgers equation associated with the initial data $v_0$. We have therefore
$$\|u(t)\|_\infty\le2\,\sqrt{\frac{2\|v_0\|_1}t\,}\,,$$
where the decay rate $t^{-\frac12}$ is independent of the space dimension. In particular this decay is faster than that given by Corollary \ref{c:dispun} when $n\ge3$.

The way this faster decay is compatible with the optimality of (\ref{eq:dispq}) is well explained by a study of the growth of the support of the solution. In the most favorable case where the data $u_0$ is bounded with compact support, the argument above yields $\|u(t)\|_\infty=O((1+t)^{-1/2})$. It is easy to infer that the width of ${\rm Supp}(u(t))$ in the $y_1$-direction expands as $O(\sqrt t\,)$ (one might have used the comparison with the solution $v$ above). Likewise, the width in the $y_2$-direction is an $O(\log t)$ and that in the other $y_k$-directions remains bounded because
$$\int_0^\infty(1+t)^{-\frac k2}dt<\infty.$$
On the contrary, if $u_0\in L^1(\R^n)$ has compact support but is not bounded by an integrable fonction $v_0(y_1)$ as above, Corollary \ref{c:dispun} gives only $\|u(t)\|_\infty=O(t^{-\kappa})$. It turns out that $n\kappa\ge1$ when $n\ge2$, and therefore 
$$\int_0t^{-n\kappa}dt=+\infty.$$
This suggest that the width of the support in the $y_n$-direction is immediately infinite: the support of $u(t)$ is unbounded for every $t>0$. The solution has a tail in the last direction, and this tail is responsible for a slow $L^\infty$-decay, at rate $t^{-\kappa}$ instead of $t^{-\frac12}$\,.

This analysis suggests in particular that the {\em fundamental solution} $U_m$, if it exists, should have an unbounded support in the space variable when $n\ge2$. The terminology denotes an entropy solution of (\ref{eq:Bru}), say a non-negative one, with the property that 
$$U_m(t)\stackrel{t\rightarrow0+}\longrightarrow m\,\delta_{y=0}$$
in the vague sense of bounded measures. In particular,
$$\int_{\R^n}U_m(t,y)\,dy\equiv m.$$
This behaviour is in strong constrast with the one-dimensional situation, where 
$$U_m(t,y)=\frac y t\,{\bf 1}_{(0,\sqrt{2m t}\,)}$$
is compactly supported at every time.

The existence of a fundamental solution is left as an open problem. It should play an important role in the time-asymptotic analysis of entropy solutions of finite mass. This asymptotics has been known in one-space dimension since the seminal works by P. Lax \cite{Lax} and C. Dafermos \cite{Daf_NW}.

\paragraph{Preliminary works.} The authors posted, separately, recent preprints on this subject in ArXiv database, see \cite{Ser_prelim,Sil_prelim}. The present paper supersedes both of them.

\paragraph{Outline of the article.}
We prove a special case of the dispersive estimate (\ref{eq:dispq}), that for the pairs $(p,p^*)$, in  Section \ref{s:pstar}. We treat the case $(p,\infty)$ in Section \ref{s:DG}. This allows us to extend the (\ref{eq:dispq}) to every pair $(p,q)$ with $p\le q$. The construction of the semi-group over every $L^p$-space is done in Section \ref{s:SGp}. We show in Section \ref{s:mono} how these ideas adapt to a scalar equation when the fluxes $f_j$ are monomials. The last section describes how the first argument, which involves Compensated Integrability, can be adapted to conservation laws with arbitrary flux.

\paragraph{Acknowledgements.} We are indebted to C. Dafermos, who led us to collaborate.

\section{Dispersive estimate ; the case $(p,p^*)$}\label{s:pstar}

To begin with, we recall that the Burgers equation enjoys an exceptional one-parameter transformation group, a fact already noted in \cite{Sil}~:
Let $u$ be an entropy solution of the Cauchy problem for (\ref{eq:Bru}) and $\lambda$ be a positive constant. Then the function
$$v(t,y)=\frac1\lambda\,u( t,\lambda y_1,\ldots,\lambda^ny_n)$$
is an entropy solution associated with the initial data
$$v_0(y)=\frac1\lambda\,u_0(\lambda y_1,\ldots,\lambda^ny_n).$$
The following identities will be used below:
\begin{eqnarray}
\label{eq:intvlat}
\int_0^\tau dt\int_{\R^n}v(t,y)^qdy & = & \lambda^{-q-\frac{d(d-1)}2}\int_0^\tau dt\int_{\R^n}u(t,y)^qdy, \\\label{eq:intvla}
\int_{\R^n}v_0(y)^qdy & = & \lambda^{-q-\frac{d(d-1)}2}\int_{\R^n}u_0(y)^qdy.
\end{eqnarray}

\bigskip

Let $u_0^\pm$ be the positive and negative parts of the initial data: $u_0^-\le u_0\le u_0^+$ with $u_0(x)\in\{u_0^-(x),u_0^+(x)\}$ everywhere. Denote $u_\pm$ the entropy solutions associated with the data $u_0^{\pm}$. By the maximum principle, we have $u_-\le u\le u_+$ everywhere. Because of $\|u(t)\|_q\le\|u_-(t)\|_q+\|u_+(t)\|_q$ and $\|u_0\|_p=(\|u_0^-\|_p^p+\|u_0^+\|_p^p)^{1/p}$, it suffices to proves the estimate for $u_\pm$, that is for initial data that are {\em signed}. And since $v(t,y)=-u(t,-y_1,y_2,\ldots,(-1)^ny_n)$ is the entropy solution associated with $v_0(y)=-u_0(-y_1,y_2,\ldots,(-1)^ny_n)$, it suffices to treat the case of a non-negative initial data.

We therefore suppose from now on that $u_0\in L^1\cap L^\infty(\R^n)$ and $u_0\ge0$, so that $u\ge0$ over $\R_+\times\R^n$. We wish to estimate $\|u(t)\|_q$ in terms of $\|u_0\|_p$ when $q=p^*=d(1+\frac pn)$. We point out that $p^*>p$.

\subsection{A Strichartz-like inequality}

If $a\in\R$, we define a symmetric matrix
$$M(a)=\left(\frac{a^{i+j+p}}{i+j+p}\right)_{0\le i,j,\le n}.$$
Remarking that
$$M(a)=\int_0^aV(s)\otimes V(s)\,s^{p-1}ds,\qquad V(s)=\begin{pmatrix} 1 \\ \vdots \\ s^n \end{pmatrix},$$
we obtain that $M(a)$ is positive definite whenever $a>0$. Obviously,
$$\det M(a)=H_{d,p}\,a^{d(p+d-1)}=H_{d,p}\,a^{np^*},$$
where 
$$H_{d,p}=\left\|\frac1{i+j+p}\right\|_{0\le i,j,\le n}>0$$
is a Hilbert-like determinant.

Let us form the symmetric tensor
$$T(t,y)=M(u(t,y)),$$
with positive semi-definite values. Its row of index $i$ is formed of $(\eta_{i+p}(u), q_{i+p}(u))$, an entropy-flux pair where $\eta_r(s)=\frac{|s|^r}r\,$ is convex. In the special case where $p=1$ and $i=0$, it is divergence-free  because of (\ref{eq:Bru}) itself.
Otherwise, it is not divergence-free in general, although it is so wherever $u$ is a classical solution.  But the entropy inequality tells us that the opposite of its divergence if a non-negative, hence bounded measure, 
$$\mu_r=-{\rm div}_{t,y}(\eta_r(u), q_r(u))\ge0.$$
The total mass of $\mu_r$ over a slab $(0,\tau)\times\R^n$ is given by
$$\|\mu_r\|=\int_{\R^n}\eta_r(u_0(y))\,dy-\int_{\R^n}\eta_r(u(\tau,y))\,dy\le\int_{\R^n}\frac{u_0(y)^r}r\,dy.$$
Since the latter bound does not depend upon $\tau$, $\mu_r$ is actually a bounded measure other $\R_+\times\R^n$.

\bigskip

We conclude that the row-wise divergence of $T$ is a (vector-valued) bounded measure, whose total mass is bounded above by
$$\sum_{j=0}^n\int_{\R^n}\frac{u_0(y)^{j+p}}{j+p}\,dy.$$
We may therefore apply Compensated Integrability (Theorems 2.2 and 2.3 of \cite{Ser_CI}) to the tensor $T$, that is
$$\int_0^\tau dt\int_{\R^n}(\det T)^{\frac1{d-1}}dy\le c_d\left(\|T_{0\bullet}(0,\cdot)\|_1+\|T_{0\bullet}(\tau,\cdot)\|_1+\|{\rm Div}_{t,y}T\|_{{\cal M}((0,\tau)\times\R^n)}\right)^{\frac d{d-1}}.$$
Because of
$$\|T_{0\bullet}(t,\cdot)\|_1=\sum_{j=0}^n\int_{\R^n}\frac{u(t,y)^{j+p}}{j+p}\,dy.\le\sum_{j=0}^n\int_{\R^n}\frac{u_0(y)^{j+p}}{j+p}\,dy.,$$
we deduce
\begin{equation}
\label{eq:nonhom}
\int_0^\tau dt\int_{\R^n}u^{p^*}dy\le c_{d,p}\left(\sum_{j=0}^n\int_{\R^n}u_0(y)^{j+p}\,dy\right)^{\frac d{d-1}}.
\end{equation}
Again, the right-hand side does not depend upon $\tau$, thus the inequality above is true also for $\tau=+\infty$.

\bigskip

The only flaw in the estimate (\ref{eq:nonhom}) is the lack of homogeneity of its right-hand side. To recover a well-balanced inequality, we use the scaling, in particular the formul\ae\, (\ref{eq:intvla}).
Applying (\ref{eq:nonhom}) to the pair $(v,v_0)$ instead, 
we get a parametrized inequality
$$\left(\int_0^\infty dt\int_{\R^n}u^{p^*}dy\right)^{\frac{d-1}d}\le c_d\lambda^{\frac{d-1}2}\sum_{j=0}^n\lambda^{-j}\int_{\R^n}u_0(y)^{j+p}\,dy,$$
where $\lambda>0$ is up to our choice.
In order to minimize the right-hand side, we select the  value
$$\lambda=\left(\int_{\R^n}u_0(y)^{n+p}dy/\int_{\R^n}u_0(y)^pdy\right)^{\frac1n}.$$
The extreme terms, for $j=0$ or $n$, contribute on a equal foot with
$$\left(\int_{\R^n}u_0(y)^{n+p}dy\right)^{\frac12}\left(\int_{\R^n}u_0(y)^pdy\right)^{\frac12}.$$
The other ones, which are
$$\left(\int_{\R^n}u_0(y)^{n+p}dy/\int_{\R^n}u_0(y)^pdy\right)^{\frac12-\frac j{d-1}}\int_{\R^n}u_0^{j+p}dy,$$
are bounded by the same quantity, because of H\"older inequality. We end therefore with the fundamental estimate of Strichartz style
\begin{equation}
\label{eq:estfond}
\left(\int_0^\infty\!\int_{\R^n}u^{p^*}dydt\right)^{\frac{d-1}d}\le c_d\left(\int_{\R^n}u_0(y)^{p+n}dy\right)^{\frac12}\left(\int_{\R^n}u_0(y)^pdy\right)^{\frac12}.
\end{equation}

\subsection{Proof of estimate (\protect\ref{eq:dispq})}\label{ss:dec}

We shall contemplate (\ref{eq:estfond}) as a differential inequality. To the end, we define
$$X(t):=\int_{\R^n}u^{p^*}dy=\|u(t)\|_{p^*}^{p^*}$$
Noticing that $p+n$ is less than $p^*$, and using H\"older inequality, we get
$$\int_{\R^n}|w|^{p+n}dy\le\left(\int_{\R^n}|w|^pdy\right)^a\left(\int_{\R^n}|w|^{p^*}dy\right)^b$$
for
$$a=\frac{p+n}{p+dn}\,,\qquad b=\frac{n^2}{p+dn}\,.$$
The inequality (\ref{eq:estfond}) implies therefore
$$\left(\int_0^\infty X(t)\,dt\right)^{\frac{2n}d}\le c_d \|u_0\|_p^{p(1+a)} X(0)^b.$$

Considering the solution $w(t,y)=u(t+\tau,y)$, whose initial data is $u(\tau,\cdot)$, we also have
\begin{equation}
\label{eq:Xtau}
\left(\int_\tau^\infty X(t)\,dt\right)^{\frac{2n}{db}}\le c_d \|u(\tau)\|_p^{p\frac{1+a}b} X(\tau)\le c_d\|u_0\|_p^{p\frac{1+a}b} X(\tau).\end{equation}
Let us denote 
$$Y(\tau):=\int_\tau^\infty X(t)\,dt.$$
We recast (\ref{eq:Xtau}) as
$$Y^\rho+c_d\|u_0\|_p^\mu Y'\le0,\qquad\rho:=\frac{2n}{db}\,\quad\mu:=p\,\frac{1+a}b\,.$$
Remark that $\rho=2\frac{p+dn}{dn}>2$.
Multiplying by $Y^{-\rho}$ and integrating, we infer
$$t+c_d\|u_0\|_p^\mu Y(0)^{1-\rho}\le c_d\|u_0\|_p^\mu Y(t)^{1-\rho}.$$
This provides a first decay estimate
$$Y(t)\le c_d\|u_0\|_p^{\frac\mu{\rho-1}}t^{-\frac1{\rho-1}}.$$
Remarking that $t\mapsto X(t)$ is a non-increasing function, so that
$$\frac\tau2\,X(\tau)\le Y(\frac\tau2),$$
we deduce the ultimate decay result
$$X(t)\le c_d\|u_0\|_p^{\frac\mu{\rho-1}}t^{-\frac\rho{\rho-1}}.$$
Restated in terms of a Lebesgue norm of $u(t)$, it says
\begin{equation}
\label{eq:decayi}
\|u(t)\|_{p^*}\le c_d\|u_0\|_p^{\alpha(p,{p^*})}\, t^{-\beta(p,{p^*})},
\end{equation}
where
$\alpha(p,q)$ and $\beta(p,q)$ are given in (\ref{eq:aldef}) and (\ref{eq:bedef}). This is a special case of (\ref{eq:dispq}).

\section{General pairs $(p,q)$ where $p<q\le\infty$}\label{s:DG}

Because of (\ref{eq:albeH}) and of the H\"older inequality, it will be enough to prove (\ref{eq:dispq}) when $q=+\infty$. Once again, it is sufficient to treat the case of non-negative data\,/\,solutions.

\subsection{An estimate for $(u-\ell)_+$}

Let $\ell>0$ be a given number. We denote $w_\ell$ the entropy solution of (\ref{eq:Bru}) associated with the initial data $(u_0-\ell)_++\ell=\max\{u_0,\ell\}$. The function $z_\ell:=w_\ell-\ell$ is an entropy solution of a modified conservation law
$$\partial_tz_\ell+\sum_{k=1}^n\partial_k\frac{(z_\ell+\ell)^{k+1}}{k+1}=0.$$
This is not exactly the Burgers equation for $z_\ell$. However the $(n+2)$-uplet $(1,X+\ell,\ldots,\ldots,\frac{(X+\ell)^{n+1}}{n+1})$ is a basis of $\R_{n+1}[X]$. We pass from this basis to $(1,X,\ldots,\frac{X^{n+1}}{n+1})$ by a triangular matrix with unit diagonal. There exists therefore a change of coordinates
$$\binom{t}{y'}=P\binom{t}{y}=\begin{pmatrix} 1 & 0 \\ \vdots & Q \end{pmatrix}\binom{t}{y},$$
where $Q$ is a unitriangular matrix, 
such that $z_\ell$ obeys the Burgers equation in the new coordinates:
$$\frac{\partial z_\ell}{\partial t}+\sum_{k=1}^n\frac{\partial}{\partial y_k'}\,\frac{z_\ell^{k+1}}{k+1}=0.$$
We may therefore apply (\ref{eq:decayi}) to $z_\ell$~:
$$\left(\int_{\R^n}z_\ell(t,y')^{p^*}dy'\right)^{\frac1{p^*}}\le c_d\left(\int_{\R^n}z_\ell(0,y')^p\,dy'\right)^{\frac{\alpha(p,{p^*})}p}\, t^{-\beta(p,{p^*})}.$$
Remarking that the time variable is unchanged, and the Jacobian of the change of variable $y\mapsto y'$ at fixed time equals one, we have actually
$$\|z_\ell(t)\|_{p^*}\le c_d\|z_\ell(0)\|_p^{\alpha(p,{p^*})}\, t^{-\beta(p,{p^*})}.$$
Finally, the maximum principle tells us that $u\le w_\ell$. The inequality above is therefore an estimate of the positive part of $u-\ell$~:
\begin{equation}
\label{eq:umsell}
\|(u-\ell)_+(t)\|_{p^*}\le c_d\|(u_0-\ell)_+\|_p^{\alpha(p,{p^*})}\, t^{-\beta(p,{p^*})}.
\end{equation}

\subsection{An iteration {\em \`a la} De Giorgi}

We now prove the $L^p$-$L^\infty$ estimate, in the special case where $\|u_0\|_p=1$. We recall that $u_0$ is non-negative.

For the moment, we fix an arbitrary constant $B>0$, which we will choose large enough in the end of the proof. Then we define the following sequences for $k\in\N$~:
$$t_k=1-2^{-k},\quad\ell_k=Bt_k,\quad w_k=(u-\ell_k)_+,\quad a_k=\|w_k(t_k)\|_p.$$
Remark that the sequences $\ell_k$ and $w_k$ are increasing and decreasing, respectively. Since $t_0=0$, we have $a_0=\|u_0\|_p=1$.

For each value of $k$, we apply (\ref{eq:umsell}) in order to estimate $\|w_{k+1}(t_{k+1})\|_{p^*}$ in terms of $\|w_{k+1}(t_k)\|_p$. For the sake of simplicity, we write $\alpha,\beta$ for $\alpha(p,p^*)$ and $\beta(p,p^*)$. We get
$$\|w_{k+1}(t_{k+1})\|_{p^*}\le c_{d,p}\|w_{k+1}(t_k)\|_p^{\alpha}\, (t_{k+1}-t_k)^{-\beta}=c_{d,p}2^{\beta(k+1)}\|w_{k+1}(t_k)\|_p^{\alpha}\le c_{d,p}2^{\beta(k+1)}a_k^{\alpha}.$$
With H\"older inequality, we have also
$$a_{k+1}=\|w_{k+1}(t_{k+1})\|_{p}\le\|w_{k+1}(t_{k+1})\left\|_{p^*}\|{\bf1}_{\{y\,:\,w_{k+1}(t_{k+1},y)>0\}}\right\|_r$$
where
$$\frac1p=\frac1{p^*}+\frac1r\,.$$
Remark that $r>1$. Combining both inequalities, we obtain
$$a_{k+1}\le c_{d,p}2^{\beta(k+1)}a_k^\alpha\left|\{y\,:\,w_{k+1}(t_{k+1},y)>0\}\right|^{\frac1r}.$$
Observing that $w_{k+1}>0$ implies $w_k>B2^{-k-1}$, we infer
$$a_{k+1}\le c_{d,p}2^{\beta(k+1)}a_k^\alpha\left|\{y\,:\,w_k(t_{k+1},y)>B2^{-k-1}\}\right|^{\frac1r}.$$
We now use Chebychev Inequality
$$\left|\{y\,:\,w_k(t_{k+1},y)>B2^{-k-1}\}\right|^{\frac1p}\le B^{-1}2^{k+1}\|w_k(t_{k+1})\|_p\le B^{-1}2^{k+1}\|w_k(t_{k})\|_p$$ 
to deduce
$$a_{k+1}\le c_{d,p}B^{-\frac pr}2^{(\beta+\frac pr)(k+1)}a_k^{\alpha+\frac pr}=C2^{Ck}a_k^{1+\delta}B^{-\gamma}.$$
We have set $\delta=\alpha-\frac p{p^*}$ and $\gamma=\frac pr$\,.

By a direct computation, we verify that $\delta$ is positive:
$$\alpha-\frac p{p^*}=\frac{p^*h(p^*)-ph(p)}{p^*h(p)}=2\,\frac{p^*-p}{p^*h(p)}>0.$$
The sequence $b_k:=B^{-\frac\gamma\delta}a_k$, which starts with $b_0=B^{-\frac\gamma\delta}$, satisfies therefore a recurrence relation
$$b_{k+1}\le C 2^{Ck}b_k^{1+\delta}.$$
It is known that if $b_0$ is small enough, that is if $B$ is large enough, then $b_k\rightarrow0+$ as $k\rightarrow+\infty$. Equivalently, $a_k\rightarrow0+$.

We have therefore found a constant $B>0$ such that
$$\|(u-\ell_k)_+(1)\|_p\le\|(u-\ell_k)_+(t_k)\|_p=a_k\rightarrow0+.$$
Since $\ell_k\rightarrow B$, this means exactly that $\|u(1)\|_\infty\le B$.

\subsection{End of the proof of dispersive estimates}

Let $u_0\in L^1\cap L^\infty(\R^n)$ be non-negative. For two positive parameters $\lambda,\mu$, the entropy solution associated with the data
$$v_0(y)=\frac1\lambda u_0(\mu\lambda y_1,\ldots,\mu\lambda^ny_n)$$
is the function
$$v(t,y)=\frac1\lambda u(\mu t,\mu\lambda y_1,\ldots,\mu\lambda^ny_n).$$
If
\begin{equation}
\label{eq:lamu}
\lambda^{p+\frac{n(n+1)}2}\mu^n=\int_{\R^n}u_0(y)^p\,dy,
\end{equation}
then $\|v_0\|_p=1$ and we may apply the previous paragraph: $\|v(1)\|_\infty\le B$. In terms of $u$, this writes
$$\|u(\mu)\|_\infty\le B\lambda.$$
Eliminating $\lambda$ with (\ref{eq:lamu}), this gives
$$\|u(\mu)\|_\infty\le B\left(\mu^{-n}\|u_0\|_p^p\right)^{\frac2{n^2+n+2p}},$$
which is nothing but the dispersive estimate (\ref{eq:dispq}) for $q=+\infty$.

\bigskip

There remains to pass from $q=+\infty$ to every $q\in[p,+\infty]$. We do that by applying the H\"older inequality. Writing
$$\frac1q=\frac{1-\theta}p+\frac\theta\infty\,,$$
we have
$$\|u(t)\|_q\le\|u(t)\|_p^{1-\theta}\|u(t)\|_\infty^\theta\le\|u_0\|_p^{1-\theta}\left(Bt^{-\beta(p,\infty)}\|u_0\|_p^{\alpha(p,\infty)}\right)^\theta.$$
We conclude by using the relations (\ref{eq:albeH}).

\section{The $L^p$-semi-group for finite exponents}\label{s:SGp}

We now prove Theorem \ref{th:wpp}. We start with a remark about $L^p$-spaces.
\begin{lemma}\label{l:app}
Let $a\in L^p(\R^n)$ be given. There exists a sequence $(b_m)_{m\ge0}$ in $(L^p\cap L^\infty)(\R^n)$,  converging towards $a$ in $L^p(\R^n)$, such that $b_m-a\in L^1(\R^n)$ and 
$$\lim_{m\rightarrow+\infty}\|b_m-a\|_1=0.$$
\end{lemma}

\bepr

Recall that 
$$L^p(\R^n)=(L^1\cap L^p)(\R^n)+(L^p\cap L^\infty)(\R^n).$$
Decomposing our function as $a=a_1+a_\infty$ where
$$a_1\in (L^1\cap L^p)(\R^n),\qquad a_\infty\in(L^p\cap L^\infty)(\R^n),$$
we may form the sequence of bounded functions $b_m:=a_\infty+\pi_m\circ a_1$, where $\pi_m$ is the projection from $\R$ onto the interval $[-m,m]$. Because of
$$\|b_m\|_p\le\|a_\infty\|_p+\|\pi_m\circ a_1\|_p\le\|a_\infty\|_p+\|a_1\|_p,$$
this sequence is bounded in $L^p(\R^n)$. In addition $b_m-a=\pi_m\circ a_1-a_1\in L^1\cap L^p(\R^n)$, and
$$\|b_m-a\|_1=\|\pi_m\circ a_1-a_1\|_1\stackrel{m\rightarrow+\infty}{\longrightarrow}0,\qquad\|b_m-a\|_p=\|\pi_m\circ a_1-a_1\|_p\stackrel{m\rightarrow+\infty}{\longrightarrow}0.$$

\enpr

\bigskip

Let $u_0\in L^p(\R^n)$ be given. In order to define $S_tu_0$, we consider a sequence $b_m$ that approximates $u_0$ in the sense of Lemma \ref{l:app}. Remark that we do not care about the construction of $b_m$, as we only use the properties stated in the Lemma.

To begin with, $u_m(t):=S_tb_m$ is well-defined and belongs to $L^\infty(\R^n)$. Because of (\ref{eq:dispq}), we have
\begin{equation}
\label{eq:umq}
\|u_m(t)\|_q\le c_{d,p,q}\|b_m\|_p^{\alpha(p,q)}\,t^{-\beta(p,q)}\le C_{p,q}(u_0)\,t^{-\beta(p,q)}.
\end{equation}
The sequence $(u_m)_{m>0}$ is thus bounded in $C_0(\tau,\infty;L^q(\R^n))$ for every $q\in[p,\infty)$ and every $\tau>0$.

The contraction property gives us
$$\|u_m(t)-u_\ell(t)\|_1\le\|b_m-b_\ell\|_1\stackrel{m,\ell\rightarrow+\infty}{\longrightarrow}0.$$
Let $r,q$ be exponents satisfying  $p\le r<q<\infty$. By H\"older inequality, we have
$$\|u_m(t)-u_\ell(t)\|_r\le\|u_m(t)-u_\ell(t)\|_1^\theta(\|u_m(t)\|_q+\|u_\ell(t)\|_q)^{1-\theta},$$
where $\theta\in(0,1]$. With (\ref{eq:umq}), we infer that 
$$\|u_m(t)-u_\ell(t)\|_r\stackrel{m,\ell\rightarrow+\infty}{\longrightarrow}0,$$
uniformly over $(\tau,\infty)$.

We have thus proved that $(u_m)_{m>0}$ is a Cauchy sequence in $C_0(\tau,\infty;L^r(\R^n))$, hence is convergent in this space. If $b_m'$ is another approximating sequence for $u_0$, and $u_m'$ the corresponding solution of the Cauchy problem, we may form an approximating sequence $c_m$ in the sense of Lemma \ref{l:app}, by alterning $b_1,b_1',b_2,b_2',\ldots$. The sequence $u_1,u_1',u_2,u_2',\ldots$ will be convergent in the sense above. This shows that the limit of $u_m$ does not depend upon the precise sequence $(b_m)_{m>0}$ chosen above. Thus we may set 
$$S_tu_0:==\lim_{m\rightarrow+\infty}u_m(t),$$
which defines a
$$u\in C_b(\R_+;L^p(\R^n))\bigcap\bigcap_{p<r<\infty}C_0(0,+\infty;L^r(\R^n)).$$

There remains to prove that $u$ is an entropy solution of (\ref{eq:Bru}). For this, we use the fact that $u_m$ is itself an entropy solution, and the convergence stated above ensures that every monomial $(u_m)^j$ in the flux $f(u_m)$, converges towards $u^j$ in $L^1_{\rm loc}$. 

The fact that $u(0)=u_0$ follows from $u_m(0)=b_m$, the $L^p$-convergence $b_m\rightarrow u_0$, and the uniform convergence $u_m(t)\rightarrow u(t)$ in $L^p(\R^n)$.

\section{Other ``monomial'' scalar conservation laws}\label{s:mono}

We consider in this section conservation laws whose fluxes are monomial. Denoting $m_k(s)=\frac{s^{k+1}}{k+1}$\,, they bear the form
\begin{equation}
\label{eq:monom}
\partial_tu+\partial_1m_{k_1}(u)+\cdots+\partial_nm_{k_n}(u)=0,
\end{equation}
where $0<k_1<\cdots<k_n$ are integers. The time derivative may be written as well $\partial_tm_{k_0}(u)$ with $k_0=0$.

As before, we may restrict to non-negative initial data $u_0$ that belong to $L^1\cap L^\infty(\R^n)$. Given an exponent $p\ge1$, our symmetric tensor if $T(t,y)=M(u(t,y))$ where now
$$M(a):=\left(m_{p+k_i+k_j-1}(a)\right)_{0\le i,j\le n}.$$
Notice that $M(a)$ is symmetric, and its upper-left entry is $\frac{a^p}p$\,. Because of
$$M(a)=\int_0^as^{p-1}V(s)\otimes V(s)\,ds,\qquad V(s):=\begin{pmatrix} s^{k_0} \\ \vdots \\ s^{k_n} \end{pmatrix},$$
it positive definite whenever $a>0$. We have 
$$\det M(a)=\Delta(p,\vec k)a^N,\qquad N=dp+2K,\quad K:=\sum_0^nk_i.$$

As above, the lines of $T$ are made of entropy-entropy flux pairs of the equation (\ref{eq:monom}). Its row-wise divergence is therefore a vector-valued bounded measure. Compensated integrability yields again an inequality
$$\left(\int_0^\infty\!dt\int_{\R^n}u(t,y)^Qdy\right)^{\frac nd}\le c_{d,p,\vec k}\sum_{j=0}^n\int_{\R^n}u_0(y)^{p+k_j}dy,\qquad Q:=\frac Nn\,.$$

The conservation law is invariant under the scaling
$$u\longmapsto v(t,y):=\frac1\lambda u(t,\lambda^{k_1}y_1,\ldots,\lambda^{k_n}y_n).$$
Applying the estimate above to $v$, we obtain a parametrized inequality~:
$$\left(\int_0^\infty\!dt\int_{\R^n}u(t,y)^Qdy\right)^{\frac nd}\le c_{d,p,\vec k}\lambda^{\frac Kd}\sum_{j=0}^n\lambda^{-k_j}\int_{\R^n}u_0(y)^{p+k_j}dy.$$
We now choose
$$\lambda=\left(\int_{\R^n}u_0(y)^{p+k_n}dy \, / \, \int_{\R^n}u_0(y)^{p}dy\right)^{\frac1{k_n}}$$ 
and obtain a Strichartz-like estimate:
$$\left(\int_0^\infty\!dt\int_{\R^n}u(t,y)^Qdy\right)^{\frac nd}\le c_{d,p,\vec k}\left(\int_{\R^n}u_0(y)^{p+k_n}dy \right)^\theta \left( \int_{\R^n}u_0(y)^{p}dy\right)^{1-\theta}$$
where
$$\theta:=\frac{K}{dk_n}\in(0,1).$$
Applying this calculation to the interval $(\tau,+\infty)$, and using the decay of the $L^p$-norm, we infer
\begin{equation}
\label{eq:monSt}
\left(\int_\tau^\infty\!dt\int_{\R^n}u(t,y)^Qdy\right)^{\frac nd}\le c_{d,p,\vec k}\left(\int_{\R^n}u(\tau,y)^{p+k_n}dy \right)^\theta \left( \int_{\R^n}u_0(y)^{p}dy\right)^{1-\theta}
\end{equation}
We may now continue the analysis with a Gronwall argument, provided $p+k_n\in(p,Q]$. We leave the interested reader to check the details. Our first dispersion estimate is
\begin{equation}
\label{eq:dispmon}
\|u(t)\|_Q\le c_{d,p}t^{-\beta(p)}\|u_0\|_p^{\alpha(p)},
\end{equation}
whenever $p\ge nk_n-2K$ (remark that for the Burgers equation, this restriction is harmless).

\bigskip

At this stage, it seems that we miss an argument in order to carry out the De Giorgi technique, because the conservation law satisfied by $u-\ell$ will be a different one. Whether it can be done here and for general conservation laws is left for a future work. What we can do at least is to combine the estimates (\ref{eq:dispmon}) in order to cover pairs $(p,q)$ of finite exponents. For instance, starting from a pair $(p,Q)$ as above and chosing $p_1=Q$, we have a corresponding $Q_1$ such that (\ref{eq:dispmon}) applies with $(p_1,Q_1)$ instead of $(p,Q)$. We infer
$$\|u(t)\|_{Q_1}\le c_{d,Q}(t/2)^{-\beta(Q)}\|u(t/2)\|_Q^{\alpha(Q)}\le c_{d,p}t^{-\beta(Q)-\alpha(Q)\beta(p)}\|u_0\|_p^{\alpha(p)\alpha(Q)}.$$
Because the iteration $p\rightarrow Q$ defines a sequence which tends to $+\infty$, and using the H\"older inequality to fill the gaps, we deduce the dispersion inequalities for the monomial conservation law:
\begin{thm}
\label{th:dispmono}
For the scalar conservation law (\ref{eq:monom}) with monomial fluxes, there exist finite constants $c_{d,p,q}$ such that whenever $p\ge nk_n-2K$, $q\in[p,\infty)$ and $u_0\in L^p\cap L^\infty(\R^n)$, we have
$$\|u(t)\|_q\le c_{d,p,q}t^{-\beta(p,q)}\|u_0\|_p^{\alpha(p,q)}.$$
The exponents are given by the formula
$$\alpha(p,q)=\frac{h(q)}{h(p)}\,,\qquad h(p):=1+\frac Kp\qquad\hbox{and}\qquad\beta(p,q)=n\left(\frac{\alpha(p,q)}p-\frac1q\right).$$
\end{thm}

As in the case of the Burgers equation, we can use these estimates in order to define the semi-group over $L^p$-spaces:
\begin{cor}\label{c:wppmono}
The semi-group $(S_t)_{t\ge0}$ for equation (\ref{eq:monom}) extends by continuity as a continuous semi-group over $L^p(\R^n)$ for every $p\in[1,+\infty)$ such taht $p\ge nk_n-2K$. It maps $L^p(\R^n)$ into $L^q(\R^n)$ for every $q\in[p,\infty)$. If $u_0\in L^p(\R^n)$, then the function $u(t,y):=(S_tu_0)(y)$ is an entropy solution with initial data $u_0$.
\end{cor}

\section{Compensated integrability for general fluxes $f$}\label{s:gen}

We consider now a multi-dimensional conservation law of the most general form (\ref{eq:mdcl}). Following the ideas developped in the Burgers and monomial cases, we begin by considering a signed, bounded initial data: $u_0\in L^1\cap L^\infty(\R^n)$, $u_0\ge0$. If $a\in\R_+$, we define a symmetric matrix
$$M_g(a)=\int_0^ag(s)Z'(s)\otimes Z'(s)\,ds,$$
where $Z(s)=(f_0(s)=s,f_1(s),\ldots,f_n(s))$ and $g$ is some positive function. This matrix is positive definite under the non-degeneracy condition that $Z([0,a])$ is not contained in an affine hyperplane. We denote
$$\Delta_g(a):=(\det M_g(a))^{\frac1n}\ge0.$$

Let us define $T(t,y):=M_\phi(u(t,y))$. Because of $u\in L^\infty(\R_+;L^1\cap L^\infty(\R^n))$, the tensor $T$ is integrable over $(0,\tau)\times\R^n$.
Each row of $T$ is made of entropy-entropy flux pairs $(F_i,Q_i)$. Since $F_i$ might not be convex, we cannot estimate the measure $\mu_i=-\partial_tF_i(u)-{\rm div}_yQ_i(u)$ directly by the integral of $F_i(u_0)$. To overcome this difficulty, we define a convex function $\phi_g$ over $\R_+$ by
$$\phi_g(0)=\phi_g'(0)=0,\qquad\phi_g''(s)=|F''(s)|,$$
where $F=(F_0,\ldots,F_n)$.
Remark that $|F'|\le\phi_g'$ and $|F|\le\phi_g$. Let $\Phi_g$ be the entropy flux associated with the entropy $\phi_g$. Then the measure $\nu_g:=-\partial_t\phi_g(u)-{\rm div}_y\Phi_g(u)$ is non-negative and a bound of its total mass is  as usual
$$\|\nu_g\|\le\int_{\R^n}\phi_g(u_0(y))\,dy.$$

We now use the kinetic formulation of (\ref{eq:mdcl}), a notion for which we refer to \cite{Per}, Theorem 3.2.1. Recall the definition of the kinetic function $\chi(\xi;a)$, whose value is ${\rm sgn}\,a$ if $\xi$ lies between $0$ and $a$, and is $0$ otherwise. There exists a non-negative bounded measure $m(t,y,\xi)$ such that the function $w(t,y,\xi)=\chi(\xi;u(t,y))$ satisfies
$$\partial_tw+f'(\xi)\cdot\nabla_yw=\frac{\partial}{\partial\xi}\,m,\qquad w(0,y;\xi)=\chi(\xi;u_0(y)).$$
If $(\eta,q)$ is an entropy-entropy flux pair, then the measure $\mu=-\partial_t\eta-{\rm div}_yq$ is given by
$$\mu=\int_\R \eta''(\xi)dm(\xi).$$
We deduce that the vector-valued measure $\mu=(\mu_0,\ldots,\mu_n)$ satisfies $|\mu|\le\nu_g$. This yields the estimate
$$\|\mu\|\le\int_{\R^n}\phi_g(u_0(y))\,dy.$$
We may therefore apply the compensated integrability, which gives here
$$\int_0^\tau dt\int_{\R^n}\Delta_g(u(t,y))\,dy\le c_d\left(\|F(u_0)\|_1+\|F(u(\tau))\|_1+\int_{\R^n}\phi_g(u_0(y))\,dy\right)^{1+\frac1n}.$$
Because of $|F|\le\phi_g$ and $\|\phi_g(u(\tau))\|_1\le\|\phi_g(u_0)\|_1$, we end up with an analog of (\ref{eq:estfond})
\begin{equation}
\label{eq:CIgen}
\int_0^\infty dt\int_{\R^n}\Delta_g(u(t,y))dy\le c_d\|\phi_g(u_0)\|_1^{1+\frac1n}.
\end{equation}

\bigskip

Whether (\ref{eq:CIgen}) can be used to prove dispersive estimates depends of the amount of nonlinearity of the equation (\ref{eq:mdcl}). We leave this question for a future work.

\end{document}